\documentclass[11pt]{article}
\usepackage{amssymb}
\usepackage{amsmath}
\usepackage{amsthm}

\newtheorem{Theorem}{Theorem}
\newtheorem{Proposition}{Proposition}
\newtheorem{Definition}{Definition}
\newtheorem{Corollary}{Corollary}

\newcommand{\pr}{\mbox{pr}_1\,}
\newcommand{\spr}{\mbox{\scriptsize pr}_1\,}
\newcommand{\op}[1]{\mathop{\oplus}\limits_{\phantom{.}#1}}
\newcommand{\opp}[2]{\mathop{\oplus}\limits_{\phantom{.}#1}^{\phantom{.}#2}}
\newcommand{\eq}[1]{\stackrel{#1}{=}}

\title{REPRESENTATIONS OF MENGER $(2,n)$-SEMIGROUPS BY MULTIPLACE FUNCTIONS}

\author{\bf Wies{\l}aw A. Dudek$^1$
and Valentin S. Trokhimenko$^2$}
\begin{document}
\sloppy \maketitle
\begin{center}
$^1$Institute of Mathematics and Computer Science,\\ Wroclaw
University of Technology,\\
50-370 Wroc{\l}aw, Poland\\
E-mail: dudek@im.pwr.wroc.pl\\
$^2$Department of Mathematics, Pedagogical University,\\
21100 Vinnitsa, Ukraine\\
E-mail: vtrokhim@sovamua.com
\end{center}

\begin{abstract}\noindent
Investigation of partial multiplace functions by algebraic methods
plays an important role in modern mathematics were we consider
various operations on sets of functions, which are naturally
defined. The basic operation for $n$-place functions is an
$(n+1)$-ary superposition $[\;\; ]$, but there are some other
naturally defined operations, which are also worth of
consideration. In this paper we consider binary Mann's
compositions $\op{1},\ldots,\op{n}$ for partial $n$-place
functions, which have many important applications for the study of
binary and $n$-ary operations. We present methods of
representations of such algebras by $n$-place functions and find
an abstract characterization of the set of $n$-place functions
closed with respect to the set-theoretic inclusion.
\end{abstract}

\textbf{1.} Let $A^n$ be the $n$-th Cartesian product of a set
$A$. The set of all partial mappings from $A^n$ into $A$ is
denoted by ${\mathcal F}(A^n,A)$. On ${\mathcal F}(A^n,A)$ we
define one $(n+1)$-ary superposition $[\phantom{=}]$ and $n$
binary compositions $\op{1},\ldots,\op{n}$ putting
\begin{eqnarray}\label{1}
&&[f,g_1,\ldots ,g_n](a_1,\ldots,a_n) =
f(g_1(a_1,\ldots,a_n),\ldots,g_n(a_1,\ldots,a_n)),\\[4pt]
\label{2}
&&(f\op{i}g)(a_1,\ldots,a_n)=f(a_1,\ldots,a_{i-1},g(a_1,\ldots,a_n),a_{i+1},\ldots,a_n),
\end{eqnarray}
where $\,a_1,\ldots,a_n\in A$, $\;f,g,g_1,\ldots,g_n\in {\mathcal
F}(A^n,A)$. It is assumed that the left and right hand of
(\ref{1}) and (\ref{2}) are defined or not defined simultaneously.
Since, as it is not difficult to verify, each composition $\op{i}$
is an associative operation, algebras of the form
$(\Phi,[\phantom{=}],\op{1},\ldots,\op{n})$, where
$\Phi\subset\mathcal{F}(A^n,A)$, are called {\it Menger
$(2,n)$-semigroups of $n$-place functions}. If $\Phi$ contains
only full $n$-place functions (called also $n$-ary
operations\footnote{In the case of universal algebras the term
"partial functions" is replaced by the term "partial operations".
Instead of "full functions" we say also "operations".}), i.e.
functions which are defined for every $(a_1,\ldots,a_n)\in A^n$,
then $(\Phi,[\phantom{=}],\op{1},\ldots,\op{n})$ is called a {\it
Menger $(2,n)$-semigroup of full $n$-place functions {\rm (or}
$n$-ary operations$)$}. An abstract characterization of such
$(2,n)$-semigroups is given in \cite{Sok}.\ $(2,2)$-semigroups of
binary operations are characterized in \cite{Tr2} and \cite{Zar}.
Algebras $(\Phi,\op{1},\ldots,\op{n})$ are considered also in
\cite{Bel} and \cite{Ya}. The algebraic properties of the
compositions $\op{1},\ldots,\op{n}$ firstly were studied by Mann
\cite{Man}.

In this paper we find an abstract characterization of the class
algebras of multiplace functions of the form $(\Phi,[\
],\op{1},\ldots,\op{n},\subset\,)$, where $\Phi$ is a set of
partial $n$-place functions and $\subset$ is an inclusion of
functions.\footnote{$f\subset g$ means that the domain of $f$ is
contained in the domain $g$ and $f(x)=g(x)$ for all $x$ from the
domain of $f$.} Such problems were earlier considered for the
class of semigroups of transformations in \cite{Sch1}, \cite{Sch2}
and for Menger algebras of multiplace functions in \cite{ST},
\cite{Tr1}.

In the theory of such algebras an important role play so-called
{\em projectors}, i.e. maps $I_i:A^n\to A$ such that
$I_i(a_1,\ldots,a_n)=a_i$ for all $a_1,\ldots,a_n\in A$ and
$i=1,\ldots,n$.\\

\textbf{2.} After Menger \cite{Meng}, an $(n+1)$-ary operation
$[\phantom{=}]$ defined on $G$ is {\em superassociative}, if for
all $x,y_1,\ldots,y_n,z_1,\ldots,z_n\in G$ the following identity
holds:
\begin{equation}\label{3}
[\,[x,y_1,\ldots,y_n],z_1,\ldots,z_n] = [x,[y_1,z_1,\ldots,z_n ]
,\ldots, [y_n,z_1,\ldots,z_n ]\,] .
\end{equation}

An $(n+1)$-ary groupoid $(G;[\phantom{=}])$ satisfying the above
identity is called a {\em Menger algebra of rank} \ $n$\ (cf.
\cite{ST}). For $n=1$ it is an arbitrary semigroup.

Such algebras are investigated by many authors. For example,
Gluskin in \cite{Gl} find some algebraic properties of such
algebras. Menger algebras satisfying some solvability criteria are
described in \cite{D1, D2, D3, ST}. Representations of such
algebras by some $n$-place functions are given in \cite{DT, ST,
Tr3}.

According to the general convention used in the theory of $n$-ary
systems, the sequence $\,x_i,x_{i+1},\ldots,x_j$, where
$i\leqslant j$, \ can be written as $\,x_i^j$ \ (for \ $i>j$ \ it
is the empty symbol). Then the above identity has the form:
\begin{equation}\label{3a}
[\,[x,y_1^n],z_1^n] = [x,[y_1,z_1^n ] ,\ldots, [y_n,z_1^n ]\,] .
\end{equation}
In this convention (\ref{1}) and (\ref{2}) can be written as
\begin{eqnarray}
&&[f,g_1^n](a_1^n) = f(g_1(a_1^n),\ldots,g_n(a_1^n)),\nonumber\\[4pt]
&&(f\op{i}g)(a_1^n)=f(a_1^{i-1},g(a_1^n),a_{i+1}^n).\nonumber
\end{eqnarray}
Similarly,\ $(\cdots((x\op{i_1}y_1)\op{i_2}y_2)\cdots)\op{i_k}y_k$
\ will be written in the abbreviated form as \
$x\op{i_1}y_1\op{i_2}\cdots\op{i_k}y_k$ \ or \
$x\opp{i_1}{i_k}y_1^k$.

Let $\{\op{1},\ldots,\op{n}\}$ be a collection of associative
binary operations defined on $G$. According to \cite{Sok} and
\cite{Ya}, $(G,\op{1},\ldots,\op{n})$ will be called a {\em
$(2,n)$-semigroup}. A $(2,n)$-semigroup with an additional
$(n+1)$-ary operation $[ \ \ ]$ satisfying (\ref{3a}) will be
called a {\em Menger $(2,n)$-semigroup} and will be denoted by
$(G,[\phantom{=}],\op{1},\ldots,\op{n})$.

Any homomorphism $P$ of a Menger $(2,n)$-semigroup $(G,[ \ \
],\op{1},\ldots,\op{n})$ onto some Menger $(2,n)$-semi\-group
$(\Phi,[ \ \ ],\op{1},\ldots,\op{n})$ of $n$-place functions, i.e.
a bijection between $G$ and $\Phi$ such that the equations
  \begin{eqnarray*}
 && P([x,y_1,\ldots,y_n])=[P(x),P(y_1),\ldots,P(y_n)],\\[4pt]
 && P(x\op{i}y)=P(x)\op{i}P(y)
  \end{eqnarray*}
are valid for all $x,y,y_1,\ldots,y_n\in G$ and $i=1,\ldots,n$,
will be called a {\em representation of $(G,[ \ \
],\op{1},\ldots,\op{n})$ by $n$-place functions}. If $P$ is an
isomorphism, then we say that $P$ is a {\em faithful
representation}.

\begin{Definition}\label{D1}
{\rm A Menger $(2,n)$-semigroup $(G,[ \ \ ],\op{1},\ldots,\op{n})$
is called {\em unitary} if it contains {\em selectors}, i.e.
elements $e_1,\ldots,e_n\in G$ such that
\begin{eqnarray}\label{4}
&&[x,e_1^n]=x,\\[4pt]
 \label{5}
 && [e_i,x_1^n]=x_i, \ \ i=1,\ldots,n,\\[4pt]
 \label{6}
&&[x,e_1^{i-1},y,e_{i+1}^n]=x\op{i}y, \ \ i=1,\ldots,n,
\end{eqnarray}
for all \ $x,x_1,\ldots,x_n,y\in G$.}
\end{Definition}

\begin{Theorem}\label{T1}
Any unitary Menger $(2,n)$-semigroup is isomorphic to some set of
full $n$-place functions containing projectors \ $I_1,\ldots,I_n$.
\end{Theorem}
\begin{proof} Let $(G,[ \ \ ],\op{1},\ldots,\op{n})$ be an unitary
Menger $(2,n)$-semigroup with selectors $e_1,\ldots,e_n$. Consider
the set $\Phi$ of $n$-ary projections $I_1,\ldots,I_n$ and
functions $\lambda_g$ defined on $G$ by the formula:
\begin{equation}\label{7}
  \lambda_g(x_1^n)=[g,x_1^n]\, .
\end{equation}
We prove that \ $P:G\to\Phi$ \ such that\ $P(g)=\lambda_g$\ is an
isomorphism.

Indeed, for all\ $g,g_1,\ldots,g_n\in G$\ and\ $x_1,\ldots,x_n\in
G$\ we have:
\begin{eqnarray}
&&\lambda_{[g,g_1^n]}(x_1^n)\eq{(\ref{7})}[[g,g_1^n],x_1^n]\eq{(\ref{3a})}[g,[g_1,x_1^n],\ldots,[g_n,x_1^n]]=
\nonumber\\[4pt]
&&\eq{(\ref{7})}\lambda_g\Big([g_1,x_1^n],\ldots,[g_n,x_1^n]\Big)\eq{(\ref{7})}
\lambda_g\Big(\lambda_{g_1}(x_1^n),\ldots,\lambda_{g_n}(x_1^n)\Big)=
\nonumber\\[4pt]
&&\eq{(\ref{1})}[\lambda_g,\lambda_{g_1},\ldots,\lambda_{g_n}](x_1^n),\nonumber
\end{eqnarray}
which implies\
$\lambda_{[g,g_1^n]}=[\lambda_g,\lambda_{g_1},\ldots,\lambda_{g_n}]$.
This proves that\ $P([g,g_1^n])= [P(g),P(g_1),\ldots,P(g_n)] $.

Similarly, for\ $g_1,g_2\in G$,\ $x_1,\ldots,x_n\in G$\ and\
$i=1,\ldots,n$\ we obtain:
\begin{eqnarray}
&&\lambda_{g_1\op{i}g_2}(x_1^n)\eq{(\ref{7})}[g_1\op{i}g_2,x_1^n]\eq{(\ref{6})}
[[g_1,e_1^{i-1},g_2,e_{i+1}^{n}],x_1^n]=\nonumber\\[4pt]
&&\eq{(\ref{3})}[g_1,[e_1,x_1^n],\ldots,[e_{i-1},x_1^n],[g_2,x_1^n],[e_{i+1},x_1^n],\ldots,
[e_n,x_1^n]]=\nonumber\\[4pt]
&&\eq{(\ref{5})}[g_1,x_1^{i-1},[g_2,x_1^n],x_{i+1}^n]\eq{(\ref{7})}
\lambda_{g_1}\Big(x_1^{i-1},\lambda_{g_2}(x_1^n),x_{i+1}^n\Big)=\nonumber\\[4pt]
&&\eq{(\ref{2})}\lambda_{g_1}\op{i}\lambda_{g_2}(x_1^n).\nonumber
\end{eqnarray}
Thus \ $\lambda_{g_1\op{i}g_2}=\lambda_{g_1}\op{i}\lambda_{g_2}$,
\ i.e. \ $P(g_1\op{i}g_2)=P(g_1)\op{i}P(g_2)$.

Moreover, $P(e_i)=I_i$ for $i=1,\dots,n$, where $e_1,\ldots,e_n$\
are selector on $G$. Indeed,
$$
\lambda_{e_i}(x_1^n)\eq{(\ref{7})}[e_i,x_1^n]\eq{(\ref{5})}x_i=I_i(x_1^n)
$$
for all $x_1,\ldots,x_n\in G$. Hence \
$P(e_i)=\lambda_{e_i}=I_i$,\ i.e.\ $\lambda_{e_i}$ is the $i$-th
projector on $G$.

Finally, if\ $\lambda_{g_1}=\lambda_{g_2}$,\ then\
$\lambda_{g_1}(e_1^n)=\lambda_{g_1}(e_1^n)$,\ consequently\
$[g_1,e_1^n]=[g_2,e_1^n]$, which (by (\ref{4})) gives\ $g_1=g_2$.
Hence\ $P:G\to\Phi$\ is an isomorphism.
\end{proof}

In the sequel, projectors of a Menger $(2,n)$-semigroup of $n$-ary
operations will be identified with selectors. For a Menger
$(2,n)$-semigroup of $n$-place functions such identification is
impossible because
$$
[I_i,g_1,\ldots,g_n] = g_i\circ\Delta_{\spr g_1\cap\ldots\cap \spr
g_n}\,
$$
where pr$_1 f$ is the domain of a function $f$ and
$\Delta_H=\{(h,h)\,|\,h\in H\}$.

Let $(G,[\phantom{=}],\op{1},\ldots,\op{n})$ be a fixed Menger
$(2,n)$-semigroup. The symbol $\mu_i(\opp{i_1}{i_s}x_1^s)$, where
$x_1,\ldots,x_s\in G$ and $\op{i_1},\ldots,\op{i_s}$ are binary
operations defined on $G$, denotes an element \
$x_{i_k}\!\opp{i_{k+1}}{i_s}\!x_{k+1}^{s}$ \ if $i=i_k$\ and\
$i\neq i_p$ for all $p<k\leqslant s$. If\ $i\neq i_p$\ for all
$i_p\in\{i_1,\ldots,i_s\}$\ this symbol is empty. For example,
$\mu_1(\op{2}x\op{1}y\op{3}z)=y\op{3}z$,
$\mu_2(\op{2}x\op{1}y\op{3}z)=x\op{1}y\op{3}z$,
$\mu_3(\op{2}x\op{1}y\op{3}z)=z$. The symbol
$\mu_4(\op{2}x\op{1}y\op{3}z)$ is empty.

\begin{Theorem}\label{T2}
A Menger $(2,n)$-semigroup \
$\mathcal{G}=(G,[\phantom{=}],\op{1},\ldots,\op{n})$ \ is
isomorphic to a Menger $(2,n)$-semigroup of $n$-place functions if
and only if it satisfies the following three identities
\begin{eqnarray}\label{8}
&&[x\op{i}y,z_1^n]=[x,z_1^{i-1},[y,z_1^n],z_{i+1}^n]\, ,\\[4pt]
 \label{9}
&&[x,y_1^n]\op{i}z=[x,y_1\op{i}z,\ldots,y_n\op{i}z]\, ,\\[4pt]
  \label{10}
&&x\opp{i_1}{i_s}y_1^s=[x,\mu_1(\opp{i_1}{i_s}y_1^s),\ldots,\mu_n(\opp{i_1}{i_s}y_1^s)]\,
,
\end{eqnarray}
where $\{i_1,\ldots,i_s\}=\{1,\ldots,n\}$, and one implication
\begin{equation}\label{11}
\bigwedge\limits_{i=1}^{n}\left(\mu_i(\opp{i_1}{i_s}x_1^s)=
\mu_i(\opp{j_1}{j_k}y_1^k)\right)\Longrightarrow
  g\opp{i_1}{i_s}x_1^s=g\opp{j_1}{j_k}y_1^k\,
\end{equation}
for all $\,g,x_1,\ldots,x_s,y_1,\ldots,y_k\in G.$
\end{Theorem}
\begin{proof} The presented proof is a modification of the proof of
Theorem 3 in \cite{Sok}, where is given the similar
characterization by $n$-ary operations.

We limit ourselves with the proof of sufficiency. Let all the
conditions of the theorem be satisfied by
$\mathcal{G}=(G,[\phantom{=}],\op{1},\ldots,\op{n})$.

 Consider the set\
$G^{*}=G\cup\{e_1,\ldots,e_n\}$,\ where\ $e_1,\ldots,e_n$\ are
different elements not belonging to $G$. For all\
$x_1,\ldots,x_s\in G$, \ $i=1,\ldots,n$, \  and operations \
$\op{i_1},\ldots,\op{i_s}$ \ defined on\ $G$\ by \
$\mu_i^{*}(\opp{i_1}{i_s}x_1^s)$ \ we denote an element of \ $G^*$
\ such that:
\begin{equation}\label{*}
  \mu_i^{*}(\opp{i_1}{i_s}x_1^s)=
\begin{cases}\phantom{==} e_i, & \mbox{if the symbol } \
\mu_i(\opp{i_1}{i_s}x_1^s) \ \mbox{ is empty}, \\[4pt]
\mu_i(\opp{i_1}{i_s}x_1^s) & \mbox{if the symbol } \
\mu_i(\opp{i_1}{i_s}x_1^s) \ \mbox{ is not empty}\, .
  \end{cases}
\end{equation}
Next, for every\ $g\in G$,\ we define on\ $G^*$\ an $n$-place
function \ $\lambda_g^*$ \ putting
$$
\lambda_g^*(x_1^n)=
  \begin{cases}
    \phantom{=}[g,x_1^n], & \mbox{if }\ x_i\in G, \ i=1,\ldots,n, \\[4pt]
    \phantom{==}g, & \mbox{if }\ x_i=e_i, \ i=1,\ldots,n,\\[4pt]
    \phantom{=}g\opp{i_1}{i_s}y_1^s, & \mbox{if }\ x_i=\mu_i^*(\opp{i_1}{i_s}y_1^s), \
    i=1,\ldots,n, \mbox{ where } y_1,\ldots,y_s\in G\, .
  \end{cases}
$$
In the other cases \ $\lambda_g^*$ \ is not defined.

Now, similarly as in \cite{Sok}, it is easy to check that \ $P:
g\mapsto\lambda_g^*$ \ is a faithful representation of \
$\mathcal{G}$.
\end{proof}

\begin{Corollary} In any unitary Menger $(2,n)$-semigroup
$(G,[\phantom{=}],\op{1},\ldots,\op{n})$ we have
\begin{equation}\label{10a}
  x\opp{i_1}{i_s}y_1^s=[x,\mu^*_1(\opp{i_1}{i_s}y_1^s),\ldots,\mu^*_n(\opp{i_1}{i_s}y_1^s)]
\end{equation}
for all $x,y_1,\ldots,y_s\in G$ and
$\,i_1,\ldots,i_s\in\{1,\ldots,n\}$.
\end{Corollary}
\begin{proof} Let
$\{j_1,\ldots,j_k\}=\{1,\ldots,n\}\setminus\{i_1,\ldots,i_s\}$.
According to the Definition \ref{D1}, in a Menger
$(2,n)$-semigroup containing selectors $e_1,\ldots,e_n$ we have\
$x\op{i}e_i=x$\ for all $x\in G$ and $i=1,\ldots,n$. Thus
$$
\mu^*_i(\opp{i_1}{i_s}y_1^s) =
\mu^*_i((\opp{i_1}{i_s}y_1^s)\op{j_1}e_{j_1}\op{j_2}\cdots\op{j_k}e_{j_k}),
$$
for $i=1,\ldots,n$.\ This implies (\ref{10}). Applying (\ref{*})
we obtain (\ref{10a}).
\end{proof}

\begin{Theorem}\label{T3}
Any Menger $(2,n)$-semigroup of $n$-place functions is isomorphic
to some Menger $(2,n)$-semigroup of full $n$-place functions.
\end{Theorem}
\begin{proof} Let
$\mathbf{F}=(\Phi,[\phantom{=}],\op{1},\ldots,\op{n})$,\ where
$\Phi\subset\mathcal{F}(A^n,A)$,\ be a Menger $(2,n)$-semigroup of
$n$-place functions. For every\ $f\in\Phi$\ we define on
$A_0=A\cup\{c\}$, where $c\not\in A$, a full $n$-place function
$f^{0}$ putting
\begin{equation}\label{12}
  f^0(x_1^n)=
  \begin{cases}
    f(x_1^n), & \mbox{if } \ x_1^n\in\pr f, \\[4pt]
    \phantom{=}c, & \mbox{if } \ x_1^n\not\in\pr f.
  \end{cases}
\end{equation}
If\  $x_1^n\in\pr[f,g_1^n]$,\ where $f,g_1,\ldots,g_n\in\Phi$,\
then\ $x_1^n\in\pr g_i$ for all $i=1,\ldots,n$\ consequently\
$(g_1(x_1^n),\ldots,g_n(x_1^n))\in\pr f$. \ Hence
\[\arraycolsep=.5mm
\begin{array}{lll}
[f,g_1^n]^0(x_1^n)&=[f,g_1^n](x_1^n)
=f(g_1(x_1^n),\ldots,g_n(x_1^n))\\[6pt]
&=f^0(g_1(x_1^n),\ldots,g_n(x_1^n))
=f^0(g^0_1(x_1^n),\ldots,g^0_n(x_1^n))\\[6pt]
&= [f^0,g^0_1,\ldots,g^0_n](x_1^n),
\end{array}
\]
which gives \
$[f,g_1^n]^0(x_1^n)=[f^0,g^0_1,\ldots,g^0_n](x_1^n)$ \ for \
$x_1^n\in\pr[f,g_1^n]$.

Now if\ $x_1^n\not\in\pr[f,g_1^n]$, then\ $[f,g_1^n]^0(x_1^n)=c$.
But $x_1^n\not\in\pr[f,g_1^n]$ implies either $x_1^n\not\in\pr
g_i$\ (for some $i=1,\ldots,n$) or
$(g_1(x_1^n),\ldots,g_n(x_1^n))\not\in\pr f$. In the first case\
$g_i^0(x_1^n)=c$, which gives
$f^0(g^0_1(x_1^n),\ldots,g^0_n(x_1^n))=c$. In the second
(according to the definition of $f^0$) we have
$f^0(g_1(x_1^n),\ldots,g_n(x_1^n))=c$, which implies
$f^0(g^0_1(x_1^n),\ldots,g^0_n(x_1^n))=c$, i.e.
$[f^0,g^0_1,\ldots,g^0_n](x_1^n)=c$. This means that
$[f,g_1^n]^0(x_1^n)=[f^0,g^0_1,\ldots,g^0_n](x_1^n)$ holds also
for\ $x_1^n\not\in\pr[f,g_1^n]$. Hence \
$[f,g_1^n]^0=[f^0,g^0_1,\ldots,g^0_n]$\ is satisfied for all\
$f,g_1,\ldots,g_n\in\Phi$.

The similar argumentation proves\ $(f\op{i}g)^0=f^0\op{i}g^0$ \
for all $f,g\in\Phi$\ and\ $i=1,\ldots,n$.

This means that the map $P:f\mapsto f^0$ is a homomorphism of
$\mathbf{F}$ onto $(\Phi_0,[\phantom{=}],\op{1},\ldots,\op{n})$,
where $\Phi_0=\{f^0\,|\,f\in\Phi\}$. In fact, as it is not
difficult to see, $P$ is an isomorphism. This completes the proof.
\end{proof}

As a simple consequence of our Theorems \ref{T2} and \ref{T3} we
obtain the main result of \cite{Sok}, which shows that a Menger
$(2,n)$-semigroup satisfying (\ref{8}), (\ref{9}), (\ref{10}) and
(\ref{11}) is isomorphic to some Menger $(2,n)$-semigroup of full
$n$-place functions.

\bigskip

\begin{Definition}
{\rm A unitary Menger $(2,n)$-semigroup
$\mathcal{G}^*=(G^{\,*},[\phantom{=}],\op{1},\ldots,\op{n})$
containing selectors $e_1,\ldots,e_n$ is called a {\em unitary
extension} of a Menger $(2,n)$-semigroup
$\,\mathcal{G}=(G,[\phantom{=}],\op{1},\ldots,\op{n})$ if:

$(a)$ \ $G\subset G^{\,*}$,

$(b)$ \ $G\cap\{e_1,\ldots,e_n\}=\emptyset$,

$(c)$ \ $G\cup\{e_1,\ldots,e_n\}$\ is the generating set of\
$\mathcal{G}^*$.}
\end{Definition}
In a unitary extension of a Menger $(2,n)$-semigroup of $n$-ary
operations selectors are replaced by projectors.

\begin{Theorem}\label{T4}
Any Menger $(2,n)$-semigroup satisfying the assumptions of
Theorem~$\ref{T2}$ can be isomorphically embedded into a unitary
extension of some Menger $(2,n)$-semigroup of full $n$-place
functions.
\end{Theorem}
\begin{proof} If a Menger $(2,n)$-semigroup
$(G,[\phantom{=}],\op{1},\ldots,\op{n})$ is isomorphic to a Menger
$(2,n)$-semigroup $(\Phi,[\phantom{=}],\op{1},\ldots,\op{n})$ of
$n$-place functions defined on $A$, then, by Theorem~\ref{T3}, it
is isomorphic to a Menger $(2,n)$-semigroup
$(\Phi_0,[\phantom{=}],\op{1},\ldots,\op{n})$ of full $n$-place
functions defined on $A_0$ (see the proof of Theorem~\ref{T3}).

Let $\,I_1,\ldots,I_n\;$ be $n$-ary projectors defined on $A_0$.
Consider the family $\left(F_k(\Phi_0)\right)_{k=0,1,\ldots}$ of
sets $F_k(\Phi_0)$ such that:

\vspace{4pt} 1. \
$F_0(\Phi_0)=\Phi_0\cup\{I_1,\ldots,I_n\}$,\vspace{2pt}

2. \ $f,g_1,\ldots,g_n\in F_k(\Phi_0)\ \Longrightarrow\
[f,g_1^n]\in F_{k+1}(\Phi_0)$,\vspace{2pt}

3. \ $f,g\in F_k(\Phi_0)\ \Longrightarrow\ f\op{i}g\in
F_{k+1}(\Phi_0 )$.
\\

It is clear that $\{ I_1,\ldots,I_n \}\cap\Phi_0 =\emptyset $ and
$\{ I_1,\ldots,I_n \}\subset F_k(\Phi_0 )$ for every $k=0,1,\ldots
$ Moreover, if $f\in F_k(\Phi_0 )$, then $ f=[ f,I_1,\ldots,I_n
]\in F_{k+1}(\Phi_0)$ (by 2). Thus $F_k(\Phi_0 )\subset
F_{k+1}(\Phi_0 )$ for every $k=0,1,\ldots$

Let\ $\Phi^*=\bigcup\limits_{k=0}^{\infty}F_k(\Phi_0 )$. For all
$f,g_1,\ldots,g_n\in\Phi^*$, there are $m_0,m_1,\ldots,m_n$ such
that $f\in F_{m_0}(\Phi_0 )$, $g_i\in F_{m_i}(\Phi_0 )$. Thus
$f,g_1,\ldots,g_n\in F_k(\Phi_0 )$ for
$k=\max\{m_0,m_1,\ldots,m_n\}$, consequently $[f,g_1^n]\in
F_{k+1}(\Phi_0 )\subset\Phi^*$. This proves that $\Phi^*$ is
closed with respect to the operation $[\phantom{=}]$.

Analogously we can verify that $\Phi^*$ is closed with respect to
$\op{1},\ldots,\op{n}$. Hence
$(\Phi^*,[\phantom{=}],\op{1},\ldots,\op{n})$ is a unitary Menger
$(2,n)$-semigroup generated by $\Phi_0\cup\{I_1,\ldots,I_n\}$. By
(\ref{12}) we have also $\Phi_0\cap\{I_1,\ldots,I_n\}=\emptyset$.
Therefore $(\Phi^*,[\phantom{=}],\op{1},\ldots,\op{n})$ is a
unitary extension of
$(\Phi_0,[\phantom{=}],\op{1},\ldots,\op{n})$.
\end{proof}

\textbf{3.} Let
$\mathcal{G}=(G,[\phantom{=}],\op{1},\ldots,\op{n})$ be a Menger
$(2,n)$-semigroup. A binary relation $\rho\subset G\times G$ is
called
\begin{itemize}
\item \textit{$v$-regular}, if
  \begin{eqnarray}\nonumber
  &&(x,y),(x_1,y_1),\ldots,(x_n,y_n)\in \rho\Longrightarrow \left\{
  \begin{array}{l}
  ([g,x_1^n],[g,y_1^n])\in\rho,\\[4pt]
  (g\op{i}x,g\op{i}y)\in\rho\end{array}\right.
  \end{eqnarray}
  for all $g,x,y,x_i,y_i\in G$, $i=1,\ldots,n$,
\item \textit{$l$-regular}, if
  \begin{eqnarray}\nonumber
  &&(x,y)\in\rho\Longrightarrow\left\{\begin{array}{l}([x,z_1^n], [y,z_1^n])\in \rho,\\[4pt]
  (x\op{i}z, y\op{i}z)\in\rho\end{array}\right.
  \end{eqnarray}
  for all $x,y,z,z_i\in G$, $i=1,\ldots,n$,
\item \textit{stable}, if
  \begin{eqnarray}\nonumber
  &&(x,y),(x_1, y_1),\ldots,(x_n,y_n)\in\rho\Longrightarrow\left\{
  \begin{array}{l}([x,x_1^n], [y,y_1^n])\in\rho,\\[4pt]
  (x_1\op{i}x_2, y_1\op{i}y_2)\in\rho\end{array}\right.
  \end{eqnarray}
  for all $x,y,x_i,y_i\in G$, $i=1,\ldots,n$.
\end{itemize}
Note that a reflexive and transitive relation $\rho$ is stable if
and only if it is $v$-regular and $l$-regular.

A subset $W$ of $G$ is called an {\em $l$-ideal} of $\mathcal G$
if $[g,w_1^{i-1},x,w_{i+1}^n]\in W$\ and $g\op{i}x\in W$\ for all
$g,w_1,\ldots,w_n\in G$,\ $i=1,\ldots,n$\ and $x\in W$.

\begin{Definition}\label{D2}{\rm By a
{\it determining pair} of a Menger $(2,n)$-semigroup $\mathcal G$
we mean an ordered pair $({\cal E}, W)$, where ${\cal E}$ is a
symmetric and transitive binary relation defined on a unitary
extension $\mathcal G^*$ of $\mathcal G$ and $\,W$ is a subset of
$G^*$ such that:
\begin{itemize}
\item[$(i)$] $G\cup\{e_1,\ldots,e_n\}\subset\pr\mathcal{E}$,
\item[$(ii)$] $\{e_1,\ldots,e_n\}\cap W=\emptyset$,
\item[$(iii)$] $[g,\mathcal{E}\langle e_1\rangle,\ldots,\mathcal{E}\langle
  e_n\rangle]\subset\mathcal{E}\langle g\rangle\,$ for all $\,g\in
  G$, where $\,\mathcal{E}\langle g\rangle\,$ denotes the $\mathcal{E}$-class containing
  $g$,
\item[$(iv)$]
  $[g,\mathcal{E}\langle\mu^*_1(\opp{i_1}{i_s}y_1^s)\rangle,\ldots,\mathcal{E}
  \langle\mu^*_n(\opp{i_1}{i_s}y_1^s)\rangle]\subset\mathcal{E}\langle
  g\opp{i_1}{i_s}y_1^s\rangle$ for all $y_1,\ldots,y_s\in G$,
  $i_1,\ldots,i_s\in\{ 1,\ldots,n\}$,
\item[$(v)$] $\mathcal{E}\cap(\mathcal{E}(G)\times\mathcal{E}(G))$
  is $v$-regular in \ $\mathcal{G}^*$,
\item[(vi)] if \ $W\neq\emptyset$, then\ $W$ is an $\mathcal{E}$-class and\ $W\cap G$\ is an
$l$-ideal of \ $\mathcal{G}$.
\end{itemize}     }
\end{Definition}

Let $\left(H_a\right)_{a\in A}$ be a collection of
$\mathcal{E}$-classes (uniquely indexed by elements of $A$) such
that $H_a\neq W$ and
$H_a\cap(G\cup\{e_1,\ldots,e_n\})\neq\emptyset$ for all $a\in A$.
Using this collection we define the set $\,\frak{A}\subset A^n\,$
in the following way:
\begin{itemize}
\item[$(a)$] if $\,g_i\in G\,$ and $\,H_{a_i}=\mathcal{E}\langle
  g_i\rangle$ for $1\leqslant i\leqslant n$, then $\;a_1^n\in\frak{A}$,
\item[$(b)$] if $\,H_{a_i}=\mathcal{E}\langle e_i\rangle$ for $1\leqslant i\leqslant n$,
  then $\;a_1^n\in\frak{A}$,
\item[$(c)$] if
$\,H_{a_i}=\mathcal{E}\langle\mu^*(\opp{i_1}{i_s}y_1^s)\rangle$
for $\leqslant i\leqslant n\,$ and some $\,y_1^s\in G$, then
  $\,a_1^n\in\frak{A}$,
\item[$(d)$] $a_1^n\in\frak{A}\,$ if and only if $\,a_1^n\,$ is
determined by $(a)$, $(b)$ or $(c)$.
\end{itemize}
Next, for any $g\in G$ we define on $A$ an $(n+1)$-ary relation
$P_{(\mathcal{E},W)}(g)$ putting:
\begin{equation}\label{23}
(a_1^n,b)\in P_{(\mathcal{E},W)}(g)\;\Longleftrightarrow\;
a_1^n\in\frak{A}\;\wedge\; [g,H_{a_1},\ldots,H_{a_n}]\subset
H_{b}\, .
\end{equation}

\noindent{\bf Remark.} If $a_1^n\in\frak{A}$ and
$[g,H_{a_1},\ldots,H_{a_n}]\subset H_b$, then $H_b\cap
G\neq\emptyset$. Indeed, if $H_{a_i}=\mathcal{E}\langle
g_i\rangle$ for $1\leqslant i\leqslant n$, where $g_i\in G$, then
$[g,g_1^n]\in H_b\cap G$. Thus $H_b\cap G\neq\emptyset$. If
$H_{a_i}=\mathcal{E}\langle e_i\rangle$ for $1\leqslant i\leqslant
n$, then $[g,e_1^n]\in H_b$ and $[g,e_1^n]=g\in G$. Therefore
$H_b\cap G\neq\emptyset$. Now if
$H_{a_i}=\mathcal{E}\langle\mu^*_i(\opp{i_1}{i_s}y_1^s)\rangle$
for $1\leqslant i\leqslant n\,$ and some $\,y_1^s\in G$, then
$\,[g,\mu^*_1(\opp{i_1}{i_s}y_1^s),\ldots,\mu^*_n(\opp{i_1}{i_s}y_1^s)]\in
H_b$. But
$\,[g,\mu^*_1(\opp{i_1}{i_s}y_1^s),\ldots,\mu^*_n(\opp{i_1}{i_s}y_1^s)]
\stackrel{(\ref{10a})}{=} g\opp{i_1}{i_s}y_1^s\in G$, which also
gives $H_b\cap G\neq\emptyset$.
\qed\\

According to the definition of $P_{(\mathcal{E},W)}(g)$
$$
(a_1^n,b),(a_1^n,c)\in P_{(\mathcal{E},W)}(g)\Longrightarrow b=c.
$$
So $P_{(\mathcal{E},W)}(g)$ is an $n$-place function such that
\begin{equation}\label{24}
  a_1^n\in\pr P_{(\mathcal{E},W)}(g)\;\Longleftrightarrow\;
  a_1^n\in\frak{A}\; \wedge\; [g,H_{a_1},\ldots,H_{a_n}]\cap
  W=\emptyset.
\end{equation}

\begin{Proposition}\label{P1}
If $(\mathcal{E},W)$ is a determining pair of a Menger
$(2,n)$-semi\-group $(G,[\phantom{=}],\op{1},\ldots,\op{n})$
satisfying the assumptions of Theorem~$\ref{T2}$, then
\begin{eqnarray}\label{25}
&&P_{(\mathcal{E},W)}([g,g_1^n])=[P_{(\mathcal{E},W)}(g),P_{(\mathcal{E},W)}(g_1),
\ldots,P_{(\mathcal{E},W)}(g_n)]\, ,\\[4pt]
 \label{26}
&&P_{(\mathcal{E},W)}(g_1\op{i}g_2)=P_{(\mathcal{E},W)}(g_1)\op{i}
P_{(\mathcal{E},W)}(g_2)\,
\end{eqnarray}
for all \ $g,g_1,\ldots,g_n\in G$ and $i=1,\ldots,n$.
\end{Proposition}
\begin{proof} We prove only (\ref{25}). The proof of (\ref{26})
is analogous.

Let $g,g_1,\ldots,g_n\in G$. If
$(a_1^n,c)\in[P_{(\mathcal{E},W)}(g),P_{(\mathcal{E},W)}(g_1),
\ldots,P_{(\mathcal{E},W)}(g_n)]$, then there are \
$b_1,\ldots,b_n\in A$\ such that
$$
(a_1^n,b_i)\in P_{(\mathcal{E},W)}(g_i)\quad {\rm and }\quad
(b_1^n,c)\in P_{(\mathcal{E},W)}(g)
$$
for $i=1,\ldots,n$. This, by (\ref{23}), is equivalent to
\begin{eqnarray}\nonumber
&&a_1^n\in\frak{A}\;\wedge\;[g_i,H_{a_1},\ldots,H_{a_n}]\subset
H_{b_i},\;\; i=1,\ldots,n,\\[4pt]
 \nonumber
&& b_1^n\in\frak{A}\;\wedge\;[g,H_{b_1},\ldots,H_{b_n}]\subset
H_c\, .
\end{eqnarray}
Let\ $h_i\in H_{b_i}$,\ $p_i\in H_{a_i}$,\ where $i=1,\ldots,n$.
The following three case are possible: 1) all $p_i$ are in $G$; 2)
$p_i=e_i$ for all $i=1,\ldots,n$; 3)
$p_i=\mu^*_i(\opp{i_1}{i_s}y_1^s)$ for all $i=1,\ldots,n$.

1) If all $p_i$ are in $G$, then\ $[g,h_1^n]\in H_c$\ and\
$[g_i,p_1^n]\in H_{b_i}\cap G$, $i=1,\ldots,n$. Since
$\big(h_i,[g_1,p_1^n]\big)\in\mathcal{E}$ for all $i=1,\ldots,n$,
the $v$-regularity of $\mathcal{E}$ implies\
$\big([g,h_1^n],[g,[g_1,p_1^n],\ldots,[g_n,p_1^n]]\big)\in\mathcal{E},$
which gives $\big([g,h_1^n],[[g,g_1^n],p_1^n]\big)\in\mathcal{E}$.
Hence $[[g,g_1^n],p_1^n]\in H_c$ consequently
$[[g,g_1^n],H_{a_1},\ldots,H_{a_n}]\subset H_c$.

2) If\ $p_1=e_1,\ldots,p_n=e_n$, then\ $[g,h_1^n]\in H_c$, $g_i\in
H_{b_i}$ and $\big(h_i, g_i\big)\in\mathcal{E}$ for
$i=1,\ldots,n$. Thus
$\big([g,h_1^n],[g,g_1^n]\big)\in\mathcal{E}$, which implies
$[g,g_1^n]\in H_c$. Hence
$[[g,g_1^n],H_{a_1},\ldots,H_{a_n}]\subset H_c$.

3) If $p_i=\mu^*_i(\opp{i_1}{i_s}y_1^s)$ for all $i=1,\ldots,n$
and some $y_1^s\in G$, then $[g,h_1^n]\in H_c$ and
$[g_i,p_1^n]\eq{(\ref{10a})}g_i\opp{i_1}{i_s}y_1^s\in H_{b_i}$ for
$i=1,\ldots,n$. Thus
$\big([g,h_1^n],[g,g_1\opp{i_1}{i_s}y_1^s,\ldots,g_n\opp{i_1}{i_s}y_1^s]\big)\in\mathcal{E},$
which, after application of (\ref{9}), proves\
$\big([g,h_1^n],[g,g_1^n]\opp{i_1}{i_s}y_1^s\big)\in\mathcal{E}$,
i.e. $\big([g,h_1^n],[[g,g_1^n],p_1^n]\big)\in\mathcal{E}$.
Therefore, as in the previous cases,\
$[[g,g_1^n],H_{a_1},\ldots,H_{a_n}]\subset H_c$.

From the above considerations we obtain
$$
a_1^n\in\frak{A} \ \wedge \
[[g,g_1^n],H_{a_1},\ldots,H_{a_n}]\subset H_c\, .
$$
Hence \ $(a_1^n,c)\in P_{(\mathcal{E},W)}([g,g_1^n])$. This
completes the proof of the inclusion
\begin{equation}\nonumber
[P_{(\mathcal{E},W)}(g),P_{(\mathcal{E},W)}(g_1),\ldots,P_{(\mathcal{E},W)}(g_n)]\subset
P_{(\mathcal{E},W)}([g,g_1^n])\, .
\end{equation}

To prove the converse inclusion consider an arbitrary element
$(a_1^n,c)$ from $P_{(\mathcal{E},W)}([g,g_1^n])$. Then $a_1^n\in
G$ and $[[g,g_1^n],H_{a_1},\ldots,H_{a_n}]\subset H_c$.

1) If $H_{a_i}=\mathcal{E}\langle h_i\rangle$ and $h_i\in G$ for
all $i=1,\ldots,n$, then $[[g,g_1^n],h_1^n]\in H_c$, which, by
(\ref{3a}), implies $[g,[g_1,h_1^n],\ldots,[g_n,h_1^n]]\in H_c$.
Obviously $[g_i,h_1^n]\in G$ for all $i$. Moreover
$H_{b_i}=\mathcal{E}\langle[g_i,h_1^n]\rangle$ for $i=1,\ldots,n$
gives $b_1^n\in\frak{A}$ and\ $[g,H_{b_1},\ldots,H_{b_n}]\subset
H_c$. It is clear that $a_1^n\in\frak{A}$ and
$[g_i,H_{a_1},\ldots,H_{a_n}]\subset H_{b_i}$ for $i=1,\ldots,n$.
Thus $(a_1^n,b_i)\in P_{(\mathcal{E},W)}(g_i)$, $i=1,\ldots n$ and
$(b_1^n,c)\in P_{(\mathcal{E},W)}(g)$. Hence
$(a_1^n,c)\in[P_{(\mathcal{E},W)}(g),P_{(\mathcal{E},W)}(g_1),
\ldots,P_{(\mathcal{E},W)}(g_n)]$.

2) If $H_{a_i}=\mathcal{E}\langle e_i\rangle$ for all
$i=1,\ldots,n$, then $[g,g_1^n]\in H_c$. For $g_i\in H_{b_i}$,
$i=1,\ldots,n$ we have $[g_i,e_1^n]\in H_{b_i}$ and
$[g_i,H_{a_1},\ldots,H_{a_n}]\subset H_{b_i}$. But
$b_1^n\in\frak{A}$ and $[g,H_{b_1},\ldots,H_{b_n}]\subset H_{c}$.
Thus $(a_1^n,b_i)\in P_{(\mathcal{E},W)}(g_i)$ and $(b_1^n,c)\in
P_{(\mathcal{E},W)}(g)$. Therefore\
$(a_1^n,c)\in[P_{(\mathcal{E},W)}(g),P_{(\mathcal{E},W)}(g_1),
\ldots,P_{(\mathcal{E},W)}(g_n)]$.

3) If
$H_{a_i}=\mathcal{E}\langle\mu^*_i(\opp{i_1}{i_s}y_1^s)\rangle$
for all $i=1,\ldots,n$ and some $y_1^s\in G$, then
$[[g,g_1^n],\mu^*_1(\opp{i_1}{i_s}y_1^s),\ldots,\mu^*_n(\opp{i_1}{i_s}y_1^s)]\in
H_c$, which, by (\ref{10a}), is equivalent to
$[g,g_1^n]\opp{i_1}{i_s}y_1^s\in H_c$. This, by (\ref{9}), implies
$[g,g_1\opp{i_1}{i_s}y_1^s,\ldots,g_n\opp{i_1}{i_s}y_1^s]\in H_c$,
where $g_i\opp{i_1}{i_s}y_1^s\in G$. But for
$H_{b_i}=\mathcal{E}\langle g_i\opp{i_1}{i_s}y_1^s\rangle$,
$i=1,\ldots,n$, we have $b_1^n\in\frak{A}$ and
$[g,H_{b_1},\ldots,H_{b_n}]\subset H_c$. Because
$g_i\opp{i_1}{i_s}y_1^s=[g_i,\mu^*_1(\opp{i_1}{i_s}y_1^s),
\ldots,\mu^*_n(\opp{i_1}{i_s}y_1^s)]$, by (\ref{10a}), then
$[g_i,H_{a_1},\ldots,H_{a_n}]\subset H_{b_i}$ for $i=1,\ldots,n$.
Hence $(a_1^n,b_i)\in P_{(\mathcal{E},W)}(g_i)$, $i=1,\ldots,n$,\
$(b_1^n,c)\in P_{(\mathcal{E},W)}(g)$ consequently
$(a_1^n,c)\in[P_{(\mathcal{E},W)}(g),P_{(\mathcal{E},W)}(g_1),
\ldots,P_{(\mathcal{E},W)}(g_n)]$.

This completes the proof the inclusion
\begin{equation}\nonumber
P_{(\mathcal{E},W)}([g,g_1^n])\subset[P_{(\mathcal{E},W)}(g),
P_{(\mathcal{E},W)}(g_1),\ldots,P_{(\mathcal{E},W)}(g_n)]
\end{equation}
and the proof of (\ref{25}).

In a similar way we prove (\ref{26}).
\end{proof}

The mapping $P_{(\mathcal{E},W)}: g\mapsto
P_{(\mathcal{E},W)}(g)$, where $g\in G$, is, by
Proposition~\ref{P1}, a representation of a Menger
$(2,n)$-semigroup $(G,[\phantom{=}],\op{1},\ldots,\op{n})$ by
$n$-place functions. This representation will be called {\em
simplest} (cf. \cite{ST}).

\medskip

Let $(P_i)_{i\in I}$ be the family of representations of a Menger
$(2,n)$-semigroup $(G,[\phantom{=}],\op{1},\ldots,\op{n})$ by
$n$-place functions defined on sets $(A_i)_{i\in I}$,
respectively. By the {\em union} of this family we mean the map
$P:g\mapsto P(g)$, where $g\in G$, and $P(g)$ is an $n$-place
function on $A=\bigcup\limits_{i\in I}A_i$ defined by
\begin{equation}\nonumber
  P(g)=\bigcup\limits_{i\in I}P_i(g)\, .
\end{equation}
If $A_i\cap A_j=\emptyset$ for all $i,j\in I$, $i\neq j$, then $P$
is called the {\em sum} of $(P_i)_{i\in I}$ and is denoted by
$P=\sum_{i\in I}P_i$. It is not difficult to see that the sum of
representations is a representation, but the union of
representations may not be a representation.

\begin{Theorem}\label{T5}
Any representation of a Menger $(2,n)$-semigroup\ $\mathcal G$\ by
$n$-place functions is a union on some family of simplest
representations of \ $\mathcal G$.
\end{Theorem}
\begin{proof} Let $P$ be a representation of a Menger
$(2,n)$-semigroup
$\mathcal{G}=(G,[\phantom{=}],\op{1},\ldots,\op{n})$ by $n$-place
functions defined on $A$, and let $\alpha\not\in A$ be some fixed
element. For every $g\in G$ we define on $A^*=A\cup\{\alpha\}$ an
$n$-place function $P^*(g)$ putting:
$$
P^*(g)(a_1^n)=
\begin{cases}
  P(g)(a_1^n), & \mbox{if } a_1^n\in\pr P(g), \\[4pt]
  \phantom{=.}\alpha, & \mbox{if } a_1^n\not\in\pr P(g).
\end{cases}
$$
It is not difficult to see that $P^*$ is a representation of
$\mathcal{G}$ by $n$-place functions defined on $A^*$, and
$P(g)\mapsto P^*(g)$, where $g\in G$, is an isomorphism of
$(P(G),[\phantom{=}],\op{1},\ldots,\op{n})$ onto
$(P^*(G),[\phantom{=}],\op{1},\ldots,\op{n})$. Because
$G\cup\{e_1,\ldots,e_n\}$ is a generating set of a unitary
extension $\mathcal{G}^*=(G^*,[\phantom{=}],\op{1},\ldots,\op{n})$
with selectors $e_1,\ldots,e_n$, then putting $P^*(e_i)=I_i$,
$i=1,\ldots,n$, where $I_i$ is the $i$-th $n$-place projector of
$A^*$, we obtain a unique extension of $P^*$ from $\mathcal{G}$ to
$\mathcal{G}^*$.

For any $a_1^n\in A^n$ we define on $G^*$ a relation
$\Theta_{a_1^n}$ such that
$$
\big(x, y\big)\in\Theta_{a_1^n}\Longleftrightarrow
P^*(x)(a_1^n)=P^*(y)(a_1^n).
$$
It is easily to verify that it is a $v$-regular equivalence
relation such that his abstract classes have the form
$H_{b}^{a_1^n}=\{x\in G^*\,|\,(a_1^n,b)\in P^*(x)\}$. The pair
$(\mathcal{E}_{a_1^n},W_{a_1^n})$, where
$$
\mathcal{E}_{a_1^n}=\Theta_{a_1^n}\cap\left(\Theta_{a_1^n}(G\cup\{e_1,\ldots,e_n\})
\times\Theta_{a_1^n}(G\cup\{e_1,\ldots,e_n\})\right),
$$
$$
W_{a_1^n}=\{x\in G^*\,|\,P^*(x)(a_1^n)=\alpha\}.
$$
is a determining pair of \ $\mathcal G$.

We prove that
\begin{equation}\nonumber
  P(g)=\bigcup\limits_{a_1^n\in A^n}P_{a_1^n}(g)
\end{equation}
for every $g\in G$, where $(P_{a_1^n})_{a_1^n\in A^n }$ is the
family of simplest representations $P_{a_1^n}$ induced by a
determining pair $(\mathcal{E}_{a_1^n},W_{a_1^n})$.

Let $(b_1^n,c)\in P(g)$. Then $g\in H_{c}^{b_1^n}$, $[g,e_1^n]=g$
and $e_i\in H_{b_i}^{b_1^n}$ for all $i=1,\ldots,n$. The
$v$-regularity of $\Theta_{a_1^n}$ implies $b_1^n\in\frak{A}$ and
$[g,H_{b_1}^{b_1^n},\ldots,H_{b_n}^{b_1^n}]\subset H_{c}^{b_1^n}$.
Hence
$$
(b_1^n,c)\in P_{b_1^n}(g)\subset\bigcup\limits_{a_1^n\in
A^n}P_{a_1^n}(g)\, .
$$
So
$$
P(g)\subset\bigcup\limits_{a_1^n\in A^n}P_{a_1^n}(g)\, .
$$

To prove the converse inclusion let \ $(b_1^n,c)\in P_{a_1^n}(g)$
\ for some\ $a_1^n\in A^n$. Then
$$
b_1^n\in\frak{A}\; \wedge\;
[g,H_{b_1}^{a_1^n},\ldots,H_{b_n}^{a_1^n}]\subset H_{c}^{a_1^n}.
$$

1) If\ $h_i\in H_{b_i}^{a_1^n}\cap G$, then $[g,h_1^n]\in
H_{c}^{a_1^n}$, consequently $(a_1^n,c)\in
P([g,h_1^n])=[P(g),P(h_1),\ldots,P(h_n)]$. Thus, for some
$d_1^n\in A^n$ we have $(a_1^n,d_i)\in P(h_i)$, $i=1,\ldots,n$ and
$(d_1^n,c)\in P(g)$. But $(a_1^n,b_i)\in P(h_i)$ implies
$b_i=d_i$. Hence $(b_1^n,c)\in P(g)$.

2) If\ $e_i\in H_{b_i}^{a_1^n}$, then $g=[g,e_1^n]\in
H_{c}^{a_1^n}$, which gives $(a_1^n,c)\in P(g)$. But for $e_i\in
H_{b_i}^{a_1^n}$, we have $(a_1^n,b_i)\in P(e_i)=P^*(e_i)=I_i$.
Thus $b_i=I_i(a_1^n)=a_i$ for $i=1,\ldots,n$. Therefore
$(b_1^n,c)\in P(g)$.

3) If\ $\mu_i^*(\opp{i_1}{i_s}y_1^s)\in H_{b_i}^{a_1^n}$ for all
$i=1,\ldots,n$ and some $y_1^s\in G$, then
$[g,\mu_1^*(\opp{i_1}{i_s}y_1^s),\ldots,\mu_n^*(\opp{i_1}{i_s}y_1^s)]\in
H_{c}^{a_1^n}$, which proves $g\opp{i_1}{i_s}y_1^s\in
H_{c}^{a_1^n}$. Thus $(a_1^n,c)\in P(g\opp{i_1}{i_s}y_1^s)=
[P(g),P(\mu_1^*(\opp{i_1}{i_s}y_1^s)),\ldots,P(\mu_n^*(\opp{i_1}{i_s}y_1^s))]$.
But there are $d_1^n\in A^n$ such that $(a_1^n,d_i)\in
P(\mu_i^*(\opp{i_1}{i_s}y_1^s))$, $i=1,\ldots,n$, $(d_1^n,c)\in
P(g)$. This together with $(a_1^n,b_i)\in
P(\mu_i^*(\opp{i_1}{i_s}y_1^s))$, $i=1,\ldots,n$, implies
$d_i=b_i$ for all $i=1,\ldots,n$.\ Therefore\ $(b_1^n,c)\in P(g)$.

Summarizing we see that in any case $P_{a_1^n}(g)\subset P(g)$.
This proves inclusion $\bigcup\limits_{a_1^n\in
A^n}P_{a_1^n}(g)\subset P(g)$\ and completes the proof of
Theorem~\ref{T5}.
\end{proof}

Note that for semigroups the analogous theorem has been proved by
B.~M.~Schein in \cite{Sch2}.

\medskip

\textbf{4.} Let
$\mathcal{G}=(G,[\phantom{=}],\op{1},\ldots,\op{n})$ be a Menger
$(2,n)$-semigroup, $x$~-- an individual variable. By $T(G)$ we
denote the set of polynomials  over $\mathcal{G}$ such that:
\begin{itemize}
\item[$(a)\;$] $x\in T(G)$,
\item[$(b)\;$] if \ $t(x)\in T(G)$, $g,g_1,\ldots,g_n\in G$,\ then\
  $[g,g_1^{i-1},t(x),g_{i+1}^n]\in T(G)$\ for all $i=1,\ldots,n$,
\item[$(c)\;$] if \ $t(x)\in T(G)$,\ $g\in G$,\ then\ $g\op{i}t(x)\in
  T(G)$\ for all $i=1,\ldots,n$,
\item[$(d)\;$] $T(G)$ contains only elements determined (a), (b) and (c).
\end{itemize}
For a nonempty subset $H$ of a Menger $(2,n)$-semigroup
$\mathcal{G}$ we define the relation $\mathcal{E}_H$ and the set
$W_H$ putting:
\begin{equation}\label{33}
  (x,y)\in\mathcal{E}_H\Longleftrightarrow\Big(\forall\;  t\in
  T(G)\Big)\Big(t(x)\in H\longleftrightarrow t(y)\in H\Big),
\end{equation}
\begin{equation}\label{34}
  x\in W_H\Longleftrightarrow \Big(\forall\;  t\in T(G)\big)\Big(t(x)\not\in H\Big).
\end{equation}
\begin{Proposition}\label{P2}
$\mathcal{E}_H$ is a $v$-regular equivalence relation on $\mathcal
G$ and $W_H$ is an $\mathcal{E}_H$-class which is an $l$-ideal or
empty set.
\end{Proposition}
\begin{proof} The fact that $\mathcal{E}_H$ is an equivalence
relation is obvious. We prove that $\mathcal{E}_H$ is $v$-regular.

Let $(x_i,y_i) \in\mathcal{E}_H$ for all $i=1,\ldots,n$ and some
$x_i,y_i\in G$. Then, by (\ref{33}), we have
\begin{equation}\label{34a}
\Big(\forall\;  t\in
  T(G)\Big)\Big(t(x_i)\in H\longleftrightarrow t(y_i)\in H\Big),
\end{equation}
for $i=1,\ldots,n$. This, for all $u,w_1,\ldots,w_n\in G$ and
polynomials of the form
$t_1\left([u,w_1^{i-1},x,w_{i+1}^n]\right)$,\ where\ $t_1(x)\in
T(G)$,\ gives
$$
\Big(t_1\in
T(G)\Big)\Big(t_1\big([u,w_1^{i-1},x_i,w_{i+1}^n]\big)\in
H\longleftrightarrow t_1\big([u,w_1^{i-1},y_i,w_{i+1}^n]\big)\in
H\Big) .
$$
Hence
\begin{equation}\label{35}
\big([u,w_1^{i-1},x_i,w_{i+1}^n], [u,w_1^{i-1},y_i,w_{i+1}^n]\big)
\in\mathcal{E}_H
\end{equation}
for all \ $i=1,\ldots,n$. In particular, for $(x_1,y_1)\in
{\mathcal{E}_H}$ we get
\begin{equation}\label{36}
 \big([u,x_1,x_2^n],[u,y_1,x_2^n])\in\mathcal{E}_H\, .
\end{equation}
Similarly, for $(x_2,y_2)\in \mathcal{E}_H$ we have
\begin{equation}\label{37}
\big([u,y_1,x_2,x_3^n],[u,y_1^2,x_3^n]\big) \in\mathcal{E}_H\, .
\end{equation}
Since $\mathcal{E}_H$ is transitive, (\ref{36}) and (\ref{37})
imply $\big([u,x_1^2,x_3^n],
[u,y_1^2,x_3^n]\big)\in\mathcal{E}_H\,$. Continuing this procedure
we obtain $\big([u,x_1^n],[u,y_1^n]\big)\in\mathcal{E}_H$.

Now let $(x,y)\in\mathcal{E}_H$, i.e. $\Big(\forall\;  t\in
T(G)\Big)\Big(t(x)\in H\longleftrightarrow t(y)\in H\Big)$. This,
for $\,t(x)=t_1(u\op{i}x)\in T(G)$, gives
$$
\Big(\forall\;\,t_1\in T(G)\Big)\Big(t_1(u\op{i}x)\in
H\longleftrightarrow t(u\op{i}y)\in H\Big).
$$
Hence $\big(u\op{i}x, u\op{i}y\big)\in\mathcal{E}_H$, which
completes the proof of the $v$-regularity of $\mathcal{E}_H$.

Assume that $W_H\neq\emptyset$. For $x,y\in W_H$, we get
$\Big(\forall\;  t\in T(G)\Big)\Big(t(x)\not\in H\Big)$ and
$\Big(\forall\;  t\in T(G)\Big)\Big(t(y)\not\in H\Big)$. Thus
$$
\Big(\forall\;  t\in T(G)\Big)\Big(t(x)\not\in H\wedge t(y)\not\in
H\Big),
$$
and, in the consequence, $\,\Big(\forall\;  t\in
T(G)\Big)\Big(t(x)\not\in H\longleftrightarrow t(y)\not\in
H\Big)$. Therefore
$$
\Big(\forall\;  t\in T(G)\Big)\Big(t(x)\in H\longleftrightarrow
t(y)\in H\Big),
$$
i.e. $\,(x,y)\in\mathcal{E}_H$. If $x\in W_H$ and
$(x,y)\in\mathcal{E}_H$, then, as it is easy to see, $y\in W_H$.
This proves that $W_H$ is an $\mathcal{E}_H$-class.

Now, according to the definition of $W_H$, for $x\in W_H$ we have
\[
\Big(\forall\;  t\in T(G)\Big)\Big(t(x)\not\in H\Big).
\]
Hence, for every $t(x)=t_1([u,w_1^{i-1},x,w_{i+1}^n])$, $u,w_i\in
G$, $i=1,\ldots,n$
\[
t_1\big([u,w_1^{i-1},x,w_{i+1}^n]\big)\not\in H .
\]
So $[u,w_1^{i-1},x,w_{i+1}^n]\in W_H$.

Analogously one can prove that $u\op{i}x\in W_H$ for all
$i=1,\ldots,n$. Therefore, $W_H$ is an $l$-ideal of $\mathcal{G}$.
\end{proof}

\medskip

Let $P$ be a fixed representation of a Menger $(2,n)$-semigroup
$\mathcal{G}$ by $n$-place functions. Consider the relation
$\zeta_P$ defined on $G$ by:
\begin{equation}\label{38}
  (g_1,g_2)\in\zeta_P\Longleftrightarrow P(g_1)\subset P(g_2)\, .
\end{equation}
It is clear that $\zeta_P$ is a quasi-order. If $P$ is a faithful
representation, then $\zeta_P$ is an order, i.e. an anti-symmetric
quasi-order. If $P$ is a sum of the family $(P_i)_{i\in I}$ of
representations of $\mathcal{G}$, then
\begin{equation}\label{38a}
  \zeta_P=\bigcap\limits_{i\in I}\zeta_{P_i}.
\end{equation}

In the case when $P=P_{(\mathcal{E},W)}$ is a simplest
representation induced by a determining pair $(\mathcal{E},W)$
(see (\ref{23}) and Proposition~\ref{P1}\,) the  quasi-order
defined by (\ref{38}) will be denoted by
$\zeta_{(\mathcal{E},W)}$.

\begin{Proposition}\label{P3}
Let $(\mathcal{E},W)$ be a determining pair of a Menger
$(2,n)$-semi\-group $\;\mathcal{G}$ satisfying the assumption of
Theorem~$\ref{T2}$. Then $\,(g_1,g_2)\in\zeta_{(\mathcal{E},W)}$
if and only if
\begin{eqnarray}\label{39}
&&\Big(\forall\;  x_1^n\in B\Big)\Big([g_1,x_1^n]\not\in
W\longrightarrow\big([g_1,x_1^n],[g_2,x_1^n]\big)\in\mathcal{E}\Big),\\[4pt]
 \label{40}
&&\Big(\forall\;  y_1^s\in
G\Big)\Big(g_1\opp{i_1}{i_s}y_1^s\not\in W\longrightarrow
\big(g_1\opp{i_1}{i_s}y_1^s,
g_2\opp{i_1}{i_s}y_1^s\big)\in\mathcal{E}\Big)
\end{eqnarray}
for \ $i_1,\ldots,i_s\in\{1,\ldots,n\}$, $s\in\mathbb{N}$,
$B=G^n\cup\{(e_1,\ldots,e_n)\}$.
\end{Proposition}
\begin{proof} If $(g_1,g_2)\in\zeta_{(\mathcal{E},W)}$, then
$P_{(\mathcal{E},W)}(g_1)\subset P_{(\mathcal{E},W)}(g_2)$. Thus
for all $a_1^n\in\frak{A}$ we have
$$
a_1^n\in\pr P_{(\mathcal{E},W)}(g_1)\longrightarrow
P_{(\mathcal{E},W)}(g_1)(a_1^n)=P_{(\mathcal{E},W)}(g_2)(a_1^n),
$$
which, by (\ref{23}) and (\ref{24}), is equivalent to
\begin{equation}\label{41}
  [g_1,H_{a_1},\ldots,H_{a_n}]\cap W=\emptyset\longrightarrow
  [g_1,H_{a_1},\ldots,H_{a_n}]\cup[g_2,H_{a_1},\ldots,H_{a_n}]
  \subset H_c
\end{equation}
for some $c\in A$. The last implication for $x_1^n\in B$ gives
(\ref{39}).

For
$x_1^n=(\mu_1^*(\opp{i_1}{i_s}y_1^s),\ldots,\mu_n^*(\opp{i_1}{i_s}y_1^s))$,
where $y_1^s\in G$, $s\in\mathbb{N}$ and
$i_1,\ldots,i_s\in\{1,\ldots,n\}$, after application of
(\ref{10a}), from (\ref{41}) we obtain
$$
g_1\opp{i_1}{i_s}y_1^s\not\in W\longrightarrow
g_1\opp{i_1}{i_s}y_1^s\in H_c\wedge g_2\opp{i_1}{i_s}y_1^s\in
H_c\,  ,
$$
and, in the consequence (\ref{40}).

In the similar way one can prove the converse statement.
\end{proof}

\medskip

We say that a Menger $(2,n)$-semigroup
$(\Phi,[\phantom{=}],\op{1},\ldots,\op{n})$ of $n$-place functions
is {\em fundamentally ordered}\ if on $\Phi$ is defined the
relation $\,\zeta_{\Phi}\,$ such that
$$
(f,g)\in\zeta_{\Phi}\Longleftrightarrow f\subset g\, .
$$
\begin{Theorem}\label{T6}
An algebraic system
$\mathcal{G}=(G,[\phantom{=}],\op{1},\ldots,\op{n},\zeta)$, where
$[\phantom{=}]$ is an $(n+1)$-ary operation,
$\op{1},\ldots,\op{n}$ are binary operations and $\zeta\subset
G\times G$, is isomorphic to a fundamentally ordered Menger
$(2,n)$-semigroup of $n$-place functions if and only if
$(G,[\phantom{=}],\op{1},\ldots,\op{n})$ is isomorphic to some
Menger $(2,n)$-semigroup of $n$-place functions and $\zeta$ is a
stable ordering relation such that
\begin{equation}\label{42}
\big(g_1,g_2\big),\, \big(g, t_1(g_1)\big),\,\big(g,
t_2(g_2)\big)\in\zeta  \Longrightarrow \big(g,
  t_2(g_1)\big)\in\zeta
\end{equation}
for all $\,g,g_1,g_2\in G$ and $t_1,t_2\in T(G)$.\\
\end{Theorem}
\begin{proof} {\it Necessity}.\ Let\
$\mathcal{G}=(G,[\phantom{=}],\op{1},\ldots,\op{n},\zeta)$ be
isomorphic to a fundamentally ordered Menger $(2,n)$-semigroup
$(\Phi,[\phantom{=}],\op{1},\ldots,\op{n},\zeta_{\Phi})$ of
$n$-place functions. We verity (\ref{42}). The rest is a
consequence of our Theorem \ref{T2}.

Assume that $g,g_1,g_2\in G$, $t_1(x),t_2(x)\in T(G)$ correspond
to $f,f_1,f_2\in\Phi$, $t_1(x),t_2(x)\in T(\Phi)$. If
$\big(g_1,g_2\big), \big(g,t_1(g_1)\big),
\big(g,t_1(g_2)\big)\in\zeta$, then $\big(f_1,f_2\big)$, $\big(f,
t_1(f_1)\big)$, $\big(f,t_2(f_2)\big)$ are in $\zeta_{\Phi}$.
Hence $f_1=f_2\circ\bigtriangleup_{\spr f_1}$, where
$\triangle_H=\{(a,a)\,|\,a\in H\}$, because $f_1\subset f_2$. By
analogy we obtain\ $f=t_1(f_1)\circ\bigtriangleup_{\spr f}$ and
$f=t_2(f_2)\circ\bigtriangleup_{\spr f}$.\ Obviously\ $f\subset
t_1(f_1)$ implies $\triangle_{\spr f}\subset\triangle_{\spr f_1}$.
Thus
\begin{eqnarray}
&&f=f\circ\bigtriangleup_{\spr
f_1}=t_2(f_2)\circ\bigtriangleup_{\spr f}\circ\bigtriangleup_{\spr
f_1}=t_2(f_2)\circ\bigtriangleup_{\spr
f_1}\circ\bigtriangleup_{\spr f}=\nonumber\\[4pt]
&&=t_2(f_2\circ\bigtriangleup_{\spr f_1})\circ\bigtriangleup_{\spr
f}=t_2(f_1)\circ\bigtriangleup_{\spr f}.\nonumber
\end{eqnarray}
Hence\ $f\subset t_2(f_1)$,\ i.e.
$\;\big(f,t_2(f_1)\big)\in\zeta_{\Phi}$.\ This proves (\ref{42}).
\\

{\it Sufficiency}. Consider the family $(P_g)_{g\in G}$ of
representations of a Menger $(2,n)$-semigroup
$(G,[\phantom{=}],\op{1},\ldots,\op{n})$, where $P_g$ is the
simplest representation induced by the determining pair
$(\mathcal{E}^*_{\zeta\langle g\rangle},W_{\zeta\langle
g\rangle})$, where $\mathcal{E}_{\zeta\langle g\rangle}$,
$W_{\zeta\langle g\rangle}$ are defined by (\ref{33}), (\ref{34}),
$\zeta\langle g\rangle=\{x\in G\,|\,\big(g, x\big)\in\zeta\}$,
$\mathcal{E}^*_{\zeta\langle g\rangle}=\mathcal{E}_{\zeta\langle
g\rangle}\cup\{(e_1,e_1),\ldots,(e_n,e_n)\}$, $e_1,\ldots,e_n$ --
selectors of $(G^*,[\phantom{=}],\op{1},\ldots,\op{n})$.

It is clear that
\begin{equation}\label{43}
  P=\sum\limits_{g\in G}P_g.
\end{equation}
is a representation of $(G,[\phantom{=}],\op{1},\ldots,\op{n})$ by
$n$-place functions.

We must proved that\ $\zeta=\zeta_P$. Let $(g_1,g_2)\in\zeta_P$,
where $g_1,g_2\in G$. Then $\,(g_1,g_2)\in\bigcap\limits_{g\in
G}\zeta_{P_g}$,\ by\ (\ref{43}) and (\ref{38a}). This together
with Proposition~\ref{P3} gives
\begin{equation}\label{44}
 \Big(\forall\;  x_1^n\in B\Big)\Big([g_1,x_1^n]\not\in W_{\zeta\langle g\rangle}
 \longrightarrow\big([g_1,x_1^n],[g_2,x_1^n]\big)\in\mathcal{E}_{\zeta\langle g\rangle}\Big)
\end{equation}
for every $g\in G$. Replacing in (\ref{44}) all $\,x_i\,$ by
$\,e_i$, $1\leqslant i\leqslant n$, we obtain the implication
\begin{equation}\nonumber
  g_1\not\in W_{\zeta\langle g\rangle}\longrightarrow
  \big(g_1,g_2\big)\in\mathcal{E}_{\zeta\langle g\rangle}\, ,
\end{equation}
which for $g=g_1$ gives $\big(g_1,
g_2\big)\in\mathcal{E}_{\zeta\langle g_1\rangle}$, because
$g_1\not\in W_{\zeta\langle g_1\rangle}$.\ Therefore
\begin{equation}\nonumber
\big(g_1,t(g_1)\big)\in\zeta \longleftrightarrow
\big(g_1,t(g_2)\big)\in\zeta\,
\end{equation}
for all $\,t\in T(G)$.\ In particular, for\ $t(x)=x$,\ we get
\begin{equation}\nonumber
  \big(g_1, g_1\big)\in\zeta \longleftrightarrow \big(g_1, g_2)\in\zeta.
\end{equation}
Applying the reflexivity of\ $\zeta$\ we obtain \ $\big(g_1,
g_2\big)\in\zeta$. \ Hence \ $\zeta_P\subset\zeta$.

Conversely, let $(g_1,g_2)\in\zeta$ for some $g_1,g_2\in G$. At
first we prove (\ref{39}). Consider $[g_1,x_1^n]\not\in
W_{\zeta\langle g\rangle}$, where $g\in G$ and $x_1^n\in B$. Then
there exists $t_1\in T(G)$ such that $\big(g,
t_1([g_1,x_1^n])\big)\in\zeta$. But $\big(g_1,g_2\big)\in\zeta$
and the stability of $\zeta$ imply $\big(t([g_1,x_1^n]),
t([g_2,x_1^n])\big)\in\zeta$ for all $t\in T(G)$ and $x_1^n\in B$.
Thus $\big(g,t([g_1,x_1^n])\big)\in\zeta$ implies $\big(g,
t([g_2,x_1^n])\big)\in\zeta$. This, together with $\big(g_1,
g_2\big)\in\zeta$, $\big(g, t_1([g_1,x_1^n])\big)\in\zeta$ and
(\ref{42}), gives\ $\big(g,t([g_1,x_1^n])\big)\in\zeta$.
Therefore\
$\big([g_1,x_1^n],[g_2,x_1^n]\big)\in\mathcal{E}_{\zeta\langle
g\rangle}$, which proves (\ref{39}).

In the similar way we can prove (\ref{40}). This, by
Proposition~\ref{P3}, shows that $(g_1,g_2)\in\zeta_{P_g}$ for
every $g\in G$, i.e. $(g_1,g_2)\in\zeta_P$. Hence
$\zeta\subset\zeta_P$, and, in the consequence $\zeta=\zeta_P$.

Now if $P(g_1)=P(g_2)$ for some $g_1,g_2\in G$, then
$P(g_1)\subset P(g_2)$ and $P(g_2)\subset P(g_1)$. Thus $\big(g_1,
g_2\big)$, $\big(g_2, g_1\big)\in\zeta$, which by the
anti-symmetry of $\zeta$ gives $g_1=g_2$. Hence $P$ is an faithful
representation and $\mathcal{G}$ is isomorphic to a fundamentally
ordered Menger $(2,n)$-semigroup
$(\Phi,[\phantom{=}],\op{1},\ldots,\op{n},\zeta_{\Phi})$ of
$n$-place functions.
\end{proof}

\end{document}